%
%

\ifx\begin\undefined\else\global\advance\srcdepth by
1\expandafter\endinput\fi

\def\begin{}
\newcount\srcdepth
\srcdepth=1

\outer\def\bye{\global\advance\srcdepth by -1
  \ifnum\srcdepth=0
    \def\endcmd{\vfill\supereject\nopagenumbers\par\vfill\supereject\end}
  \else\def\endcmd{}\fi
  \endcmd
}



\def\initialize#1#2#3#4#5#6{
  \ifnum\srcdepth=1
  \magnification=#1
  \hsize = #2
  \vsize = #3
  \hoffset=#4
  \advance\hoffset by -\hsize
  \divide\hoffset by 2
  \advance\hoffset by -1truein
  \voffset=#5
  \advance\voffset by -\vsize
  \divide\voffset by 2
  \advance\voffset by -1truein
  \advance\voffset by #6
  \baselineskip=13pt
  \emergencystretch = 0.05\hsize
  \fi
}

\def\print{\initialize{1095}
  {5.5truein}{8.5truein}{8.5truein}{11truein}{-.0625truein}}

\newif\ifblackboardbold

\blackboardboldtrue


\font\titlefont=cmbx12 scaled\magstephalf
\font\sectionfont=cmbx12

\font\scriptit=cmti10 at 7pt
\font\scriptsl=cmsl10 at 7pt
\scriptfont\itfam=\scriptit
\scriptfont\slfam=\scriptsl


\newfam\bboldfam
\ifblackboardbold
\font\tenbbold=msbm10
\font\sevenbbold=msbm7
\font\fivebbold=msbm5
\textfont\bboldfam=\tenbbold
\scriptfont\bboldfam=\sevenbbold
\scriptscriptfont\bboldfam=\fivebbold
\def\bbold{\fam\bboldfam\tenbbold}
\else
\def\bbold{\bf}
\fi


\newfam\msamfam
\font\tenmsam=msam10
\font\sevenmsam=msam7
\font\fivemsam=msam5
\textfont\msamfam=\tenmsam
\scriptfont\msamfam=\sevenmsam
\scriptscriptfont\msamfam=\fivemsam

\newfam\msbmfam
\font\tenmsbm=msam10
\font\sevenmsbm=msam7
\font\fivemsbm=msam5
\textfont\msbmfam=\tenmsbm
\scriptfont\msbmfam=\sevenmsbm
\scriptscriptfont\msbmfam=\fivemsbm

\newcount\amsfamcount 
\newcount\classcount   
\newcount\positioncount
\newcount\codecount
\newcount\n             
\def\newsymbol#1#2#3#4#5{               
\n="#2                                  
\ifnum\n=1 \amsfamcount=\msamfam\else   
\ifnum\n=2 \amsfamcount=\msbmfam\else   
\ifnum\n=3 \amsfamcount=\eufmfam
\fi\fi\fi
\multiply\amsfamcount by "100           
\classcount="#3                 
\multiply\classcount by "1000           
\positioncount="#4#5            
\codecount=\classcount                  
\advance\codecount by \amsfamcount      
\advance\codecount by \positioncount
\mathchardef#1=\codecount}              


\font\Arm=cmr9
\font\Ai=cmmi9
\font\Asy=cmsy9
\font\Abf=cmbx9
\font\Brm=cmr7
\font\Bi=cmmi7
\font\Bsy=cmsy7
\font\Bbf=cmbx7
\font\Crm=cmr6
\font\Ci=cmmi6
\font\Csy=cmsy6
\font\Cbf=cmbx6

\ifblackboardbold
\font\Abbold=msbm10 at 9pt
\font\Bbbold=msbm7
\font\Cbbold=msbm5 at 6pt
\fi

\def\small{%
\textfont0=\Arm \scriptfont0=\Brm \scriptscriptfont0=\Crm
\textfont1=\Ai \scriptfont1=\Bi \scriptscriptfont1=\Ci
\textfont2=\Asy \scriptfont2=\Bsy \scriptscriptfont2=\Csy
\textfont\bffam=\Abf \scriptfont\bffam=\Bbf \scriptscriptfont\bffam=\Cbf
\def\rm{\fam0\Arm}\def\mit{\fam1}\def\oldstyle{\fam1\Ai}%
\def\bf{\fam\bffam\Abf}%
\ifblackboardbold
\textfont\bboldfam=\Abbold
\scriptfont\bboldfam=\Bbbold
\scriptscriptfont\bboldfam=\Cbbold
\def\bbold{\fam\bboldfam\Abbold}%
\fi
\rm
}








\newlinechar=`@
\def\forwardmsg#1#2#3{\immediate\write16{@*!*!*!* forward reference should
be: @\noexpand\forward{#1}{#2}{#3}@}}
\def\nodefmsg#1{\immediate\write16{@*!*!*!* #1 is an undefined reference@}}

\def\forwardsub#1#2{\def\newref{{#2}{#1}}}

\def\forward#1#2#3{%
\expandafter\expandafter\expandafter\forwardsub\expandafter{#3}{#2}
\expandafter\ifx\csname#1\endcsname\relax\else%
\expandafter\ifx\csname#1\endcsname\newref\else%
\forwardmsg{#1}{#2}{#3}\fi\fi%
\expandafter\let\csname#1\endcsname\newref}

\def\firstarg#1{\expandafter\argone #1}\def\argone#1#2{#1}
\def\secondarg#1{\expandafter\argtwo #1}\def\argtwo#1#2{#2}

\def\ref#1{\expandafter\ifx\csname#1\endcsname\relax
  {\nodefmsg{#1}\bf`#1'}\else
  \expandafter\firstarg\csname#1\endcsname
  ~\expandafter\secondarg\csname#1\endcsname\fi}

\def\refs#1{\expandafter\ifx\csname#1\endcsname\relax
  {\nodefmsg{#1}\bf`#1'}\else
  \expandafter\firstarg\csname #1\endcsname
  s~\expandafter\secondarg\csname#1\endcsname\fi}

\def\refn#1{\expandafter\ifx\csname#1\endcsname\relax
  {\nodefmsg{#1}\bf`#1'}\else
  \expandafter\secondarg\csname #1\endcsname\fi}



\def\widow#1{\vskip 0pt plus#1\vsize\goodbreak\vskip 0pt plus-#1\vsize}



\def\marginlabel#1{}

\def\showlabelsabove{
\font\labelfont=cmss10 at 6pt
\def\marginlabel##1{\rlap{\smash{\raise 10pt\hbox{\labelfont##1}}}}
}

\newcount\seccount
\newcount\proccount
\seccount=0
\proccount=0

\def\stdskip{\vskip 9pt plus3pt minus 3pt}
\def\stdbreak{\par\removelastskip\penalty-100\stdskip}

\def\proof{\stdbreak\noindent{\sl Proof. }}

\def\qed{\vrule height 1.2ex width .9ex depth .1ex}

\def\Box{
  \ifmmode\eqno\qed
  \else\ifvmode\removelastskip\line{\hfil\qed}
  \else\unskip\quad\hskip-\hsize
    \hbox{}\hskip\hsize minus 1em\qed\par
  \fi\stdbreak\fi}

\def\references{
  \removelastskip
  \widow{.05}
  \vskip 24pt plus 6pt minus 6 pt
  \parindent=0pt
  \frenchspacing
  \leftline{\sectionfont References}
  \nobreak\stdskip\noindent}

\def\ifempty#1#2\endB{\ifx#1\endA}
\def\makeref#1#2#3{\ifempty#1\endA\endB\else\forward{#1}{#2}{#3}\fi}

\outer\def\section#1 #2\par{
  \removelastskip
  \global\advance\seccount by 1
  \global\proccount=0\relax
                \edef\numtoks{\number\seccount}
  \makeref{#1}{Section}{\numtoks}
  \widow{.05}
  \vskip 24pt plus 6pt minus 6 pt
  \message{#2}
  \leftline{\marginlabel{#1}\sectionfont\numtoks\quad #2}
  \nobreak\stdskip}

\def\proclamation#1#2{
  \outer\expandafter\def\csname#1\endcsname##1 ##2\par{
  \stdbreak
  \global\advance\proccount by 1
  \edef\numtoks{\number\seccount.\number\proccount}
  \makeref{##1}{#2}{\numtoks}
  \noindent{\marginlabel{##1}\bf #2 \numtoks\enspace}
  {\sl##2\par}
  \stdbreak}}

\def\othernumbered#1#2{
  \outer\expandafter\def\csname#1\endcsname##1{
  \stdbreak
  \global\advance\proccount by 1
  \edef\numtoks{\number\seccount.\number\proccount}
  \makeref{##1}{#2}{\numtoks}
  \noindent{\marginlabel{##1}\bf #2 \numtoks\enspace}}}

\proclamation{definition}{Definition}
\proclamation{lemma}{Lemma}
\proclamation{proposition}{Proposition}
\proclamation{theorem}{Theorem}
\proclamation{corollary}{Corollary}
\proclamation{conjecture}{Conjecture}

\othernumbered{example}{Example}
\othernumbered{remark}{Remark}
\othernumbered{construction}{Construction}
\othernumbered{problem}{Problem}

\def\figure#1{
 \global\advance\figcount by 1
 \goodbreak
 \midinsert#1\smallskip
 \centerline{Figure~\number\figcount}
 \endinsert}

\def\capfigure#1#2{
 \global\advance\figcount by 1
 \goodbreak
 \midinsert#1\smallskip
 \vbox{\small\noindent {\bf Figure~\number\figcount:} #2}
 \endinsert}

\def\capfigurepair#1#2#3#4{
 \goodbreak
 \midinsert
 #1\smallskip
 \global\advance\figcount by 1
 \vbox{\small\noindent {\bf Figure~\number\figcount:} #2}
 \vskip 12pt
 #3\smallskip
 \global\advance\figcount by 1
 \vbox{\small\noindent {\bf Figure~\number\figcount:} #4}
 \endinsert}


\def\baretable#1#2{
\vbox{\offinterlineskip\halign{
 \strut\kern #1\hfil##\kern #1
 &&\kern #1\hfil##\kern #1\cr
 #2
}}}

\def\gridtablesub#1#2#3{
\vbox{\offinterlineskip\halign{
 \strut\vrule\kern #1\hfil##\hfil\kern #2\vrule
 &&\kern #1\hfil##\kern #2\vrule\cr
 \noalign{\hrule}
 #3
 \noalign{\hrule}
}}}






\newif\iftextures

\newcount\figcount
\figcount=0
\newcount\figxscale
\newcount\figyscale
\newcount\figxoffset
\newcount\figyoffset
\newbox\drawing
\newcount\drawbp
\newdimen\drawx
\newdimen\drawy
\newdimen\ngap
\newdimen\sgap
\newdimen\wgap
\newdimen\egap

\def\drawbox#1#2#3{\vbox{
  \epsfgetbb{#2.eps} 
  \drawbp=\epsfurx
  \advance\drawbp by-\epsfllx\relax
  \multiply\drawbp by #1
  \divide\drawbp by 100
  \drawx=\drawbp bp
  \drawbp=\epsfury
  \advance\drawbp by-\epsflly\relax
  \multiply\drawbp by #1
  \divide\drawbp by 100
  \drawy=\drawbp bp
  \iftextures
  		\figxscale=#1
    \multiply\figxscale by 10
    \setbox\drawing=\vbox to \drawy{\vfil
      \hbox to \drawx{\special{illustration #2.eps scaled
\number\figxscale}\hfil}}
  \else 
    \figxoffset=-\epsfllx
    \multiply\figxoffset by#1
    \divide\figxoffset by100
    \figyoffset=-\epsflly
    \multiply\figyoffset by#1
    \divide\figyoffset by100
    \setbox\drawing=\vbox to \drawy{\vfil
      \hbox to \drawx{\includegraphics{#2.eps}\hfil}}
  \fi
  \setbox\drawing=\vbox{\offinterlineskip\box\drawing\kern 0pt}
   \ngap=0pt \sgap=0pt \wgap=0pt \egap=0pt
  \setbox0=\vbox{\offinterlineskip
    \box\drawing \ifgridlines\drawgrid\drawx\drawy\fi #3}
  \kern\ngap\hbox{\kern\wgap\box0\kern\egap}\kern\sgap}}

\def\draw#1#2#3{
  \setbox\drawing=\drawbox{#1}{#2}{#3}
  \global\advance\figcount by 1
  \edef\numtoks{\number\figcount}
  \makeref{fig:#2}{Figure}{\numtoks}
  \goodbreak
  \midinsert
  \centerline{\ifgridlines\boxgrid\drawing\fi\box\drawing}
  \smallskip
  \vbox{\offinterlineskip
    \centerline{Figure~\number\figcount}
    \smash{\marginlabel{#2}}}
  \endinsert}

\def\capdraw#1#2#3#4{
  \setbox\drawing=\drawbox{#1}{#2}{#3}
  \global\advance\figcount by 1
  \edef\numtoks{\number\figcount}
  \makeref{fig:#2}{Figure}{\numtoks}
  \goodbreak
  \midinsert
  \centerline{\ifgridlines\boxgrid\drawing\fi\box\drawing}
  \smallskip
  \vbox{\offinterlineskip
    \vskip 4pt
    \vbox{\centerline{\lineskip=3pt\small\noindent
                       {\bf Figure~\number\figcount:} #4}}
    \smash{\marginlabel{fig:#2}}}
  \endinsert}

\def\capdrawpair#1#2#3#4#5#6#7#8{
  \goodbreak
  \midinsert
  \setbox\drawing=\drawbox{#1}{#2}{#3}
  \global\advance\figcount by 1
  \edef\numtoks{\number\figcount}
  \makeref{fig:#2}{Figure}{\numtoks}
  \centerline{\ifgridlines\boxgrid\drawing\fi\box\drawing}
  \smallskip
  \vbox{\offinterlineskip
    \vskip 4pt
    \vbox{\lineskip=3pt\small\noindent {\bf Figure~\number\figcount:} #4}
    \smash{\marginlabel{fig:#2}}}
  \vskip 12pt
  \setbox\drawing=\drawbox{#5}{#6}{#7}
  \global\advance\figcount by 1
  \edef\numtoks{\number\figcount}
  \makeref{fig:#6}{Figure}{\numtoks}
  \centerline{\ifgridlines\boxgrid\drawing\fi\box\drawing}
  \smallskip
  \vbox{\offinterlineskip
    \vskip 4pt
    \vbox{\lineskip=3pt\small\noindent {\bf Figure~\number\figcount:} #8}
    \smash{\marginlabel{fig:#6}}}
  \endinsert}

\def\nextfigtoks{%
  \advance\figcount by 1%
  \edef\numtoks{\number\figcount}%
  \advance\figcount by -1}

\newif\ifgridlines
\newbox\figtbox
\newbox\figgbox
\newdimen\figtx
\newdimen\figty

\newdimen\bwd
\bwd=2sp 

\def\hline#1{\vbox{\smash{\hbox to #1{\leaders\hrule height \bwd\hfil}}}}

\def\vline#1{\hbox to 0pt{%
  \hss\vbox to #1{\leaders\vrule width \bwd\vfil}\hss}}

\def\clap#1{\hbox to 0pt{\hss#1\hss}}
\def\vclap#1{\vbox to 0pt{\offinterlineskip\vss#1\vss}}

\def\hstutter#1#2{\hbox{%
  \setbox0=\hbox{#1}%
  \hbox to #2\wd0{\leaders\box0\hfil}}}

\def\vstutter#1#2{\vbox{
  \setbox0=\vbox{\offinterlineskip #1}
  \dp0=0pt
  \vbox to #2\ht0{\leaders\box0\vfil}}}

\def\crosshairs#1#2{
  \dimen1=.002\drawx
  \dimen2=.002\drawy
  \ifdim\dimen1<\dimen2\dimen3\dimen1\else\dimen3\dimen2\fi
  \setbox1=\vclap{\vline{2\dimen3}}
  \setbox2=\clap{\hline{2\dimen3}}
  \setbox3=\hstutter{\kern\dimen1\box1}{4}
  \setbox4=\vstutter{\kern\dimen2\box2}{4}
  \setbox1=\vclap{\vline{4\dimen3}}
  \setbox2=\clap{\hline{4\dimen3}}
  \setbox5=\clap{\copy1\hstutter{\box3\kern\dimen1\box1}{6}}
  \setbox6=\vclap{\copy2\vstutter{\box4\kern\dimen2\box2}{6}}
  \setbox1=\vbox{\offinterlineskip\box5\box6}
  \smash{\vbox to #2{\hbox to #1{\hss\box1}\vss}}}

\def\boxgrid#1{\rlap{\vbox{\offinterlineskip
  \setbox0=\hline{\wd#1}
  \setbox1=\vline{\ht#1}
  \smash{\vbox to \ht#1{\offinterlineskip\copy0\vfil\box0}}
  \smash{\vbox{\hbox to \wd#1{\copy1\hfil\box1}}}}}}

\def\drawgrid#1#2{\vbox{\offinterlineskip
  \dimen0=\drawx
  \dimen1=\drawy
  \divide\dimen0 by 10
  \divide\dimen1 by 10
  \setbox0=\hline\drawx
  \setbox1=\vline\drawy
  \smash{\vbox{\offinterlineskip
    \copy0\vstutter{\kern\dimen1\box0}{10}}}
  \smash{\hbox{\copy1\hstutter{\kern\dimen0\box1}{10}}}}}

\def\figtext#1#2#3#4#5{
  \setbox\figtbox=\vbox{\hbox{#5}\kern 0pt}
  \figtx=-#3\wd\figtbox \figty=-#4\ht\figtbox
  \advance\figtx by #1\drawx \advance\figty by #2\drawy
  \dimen0=\figtx \advance\dimen0 by\wd\figtbox \advance\dimen0 by-\drawx
  \ifdim\dimen0>\egap\global\egap=\dimen0\fi
  \dimen0=\figty \advance\dimen0 by\ht\figtbox \advance\dimen0 by-\drawy
  \ifdim\dimen0>\ngap\global\ngap=\dimen0\fi
  \dimen0=-\figtx
  \ifdim\dimen0>\wgap\global\wgap=\dimen0\fi
  \dimen0=-\figty
  \ifdim\dimen0>\sgap\global\sgap=\dimen0\fi
  \smash{\rlap{\vbox{\offinterlineskip
    \hbox{\hbox to \figtx{}\ifgridlines\boxgrid\figtbox\fi\box\figtbox}
    \vbox to \figty{}
    \ifgridlines\crosshairs{#1\drawx}{#2\drawy}\fi
    \kern 0pt}}}}


\def\hpad#1#2#3{\hbox{\kern #1\hbox{#3}\kern #2}}
\def\vpad#1#2#3{\setbox0=\hbox{#3}\vbox{\kern #1\box0\kern #2}}




\def\stack#1#2#3{\vbox{\offinterlineskip
  \setbox2=\hbox{#2}
  \setbox3=\hbox{#3}
  \dimen0=\ifdim\wd2>\wd3\wd2\else\wd3\fi
  \hbox to \dimen0{\hss\box2\hss}
  \kern #1
  \hbox to \dimen0{\hss\box3\hss}}}


\def\hexp#1{%
  \setbox0=\hbox{${}^{#1}$}%
  \hbox to .5\wd0{\box0\hss}}

\def\hsub#1{%
  \setbox0=\hbox{${}_{#1}$}%
  \hbox to .5\wd0{\box0\hss}}



\def\bmatrix#1#2{{\left(\vcenter{\halign
  {&\kern#1\hfil$##\mathstrut$\kern#1\cr#2}}\right)}}

\def\rightarrowmat#1#2#3{
  \setbox1=\hbox{\small\kern#2$\bmatrix{#1}{#3}$\kern#2}
  \,\vbox{\offinterlineskip\hbox to\wd1{\hfil\copy1\hfil}
    \kern 3pt\hbox to\wd1{\rightarrowfill}}\,}

\def\leftarrowmat#1#2#3{
  \setbox1=\hbox{\small\kern#2$\bmatrix{#1}{#3}$\kern#2}
  \,\vbox{\offinterlineskip\hbox to\wd1{\hfil\copy1\hfil}
    \kern 3pt\hbox to\wd1{\leftarrowfill}}\,}

\def\rightarrowbox#1#2{
  \setbox1=\hbox{\kern#1\hbox{\small #2}\kern#1}
  \,\vbox{\offinterlineskip\hbox to\wd1{\hfil\copy1\hfil}
    \kern 3pt\hbox to\wd1{\rightarrowfill}}\,}

\def\leftarrowbox#1#2{
  \setbox1=\hbox{\kern#1\hbox{\small #2}\kern#1}
  \,\vbox{\offinterlineskip\hbox to\wd1{\hfil\copy1\hfil}
    \kern 3pt\hbox to\wd1{\leftarrowfill}}\,}








\def\quiremacro#1#2#3#4#5#6#7#8#9{
  \expandafter\def\csname#1\endcsname##1{
  \ifnum\srcdepth=1
  \magnification=#2
  \input quire
  \hsize=#3
  \vsize=#4
  \htotal=#5
  \vtotal=#6
  \shstaplewidth=#7
  \shstaplelength=#8
  \hoffset=\htotal
  \advance\hoffset by -\hsize
  \divide\hoffset by 2
  \ifnum\vsize<\vtotal
    \voffset=\vtotal
    \advance\voffset by -\vsize
    \divide\voffset by 2
  \fi
  \advance\voffset by #9
  \shhtotal=2\htotal
  \baselineskip=13pt
  \emergencystretch = 0.05\hsize
  \horigin=0.0truein
  \vorigin=0.0truein
  \shthickness=0pt
  \shoutline=0pt
  \shcrop=0pt
  \shvoffset=-1.0truein
  \ifnum##1>0\quire{#1}\else\qtwopages\fi
  \fi
}}



\quiremacro{letterbooklet} 
{1000}{4.79452truein}{7truein}{5.5truein}{8.5truein}{0.01pt}{0.66truein}
{-.0625truein}

\quiremacro{Afourbooklet}
{1095}{5.25truein}{7truein}{421truept}{595truept}{0.01pt}{0.66truein}
{-.0625truein}

\quiremacro{legalbooklet}
{1095}{5.25truein}{7truein}{7.0truein}{8.5truein}{0.01pt}{0.66truein}
{-.0625truein}

\quiremacro{twoupsub} 
{895}{4.5truein}{7truein}{5.5truein}{8.5truein}{0pt}{0pt}{.0625truein}


\quiremacro{Afourviewsub} 
{1000}{5.0228311in}{7.7625571in}{421truept}{595truept}{0.1pt}{0.5\vtotal}
{-.0625truein}


\quiremacro{viewsub}
{1095}{5.5truein}{8.5truein}{461truept}{666truept}{0.1pt}{0.5\vtotal}
{-.125truein}


\newcount\countA
\newcount\countB
\newcount\countC

\def\monthname{\begingroup
  \ifcase\number\month
    \or January\or February\or March\or April\or May\or June\or
    July\or August\or September\or October\or November\or December\fi
\endgroup}

\def\dayname{\begingroup
  \countA=\number\day
  \countB=\number\year
  \advance\countA by 0 
  \advance\countA by \ifcase\month\or
    0\or 31\or 59\or 90\or 120\or 151\or
    181\or 212\or 243\or 273\or 304\or 334\fi
  \advance\countB by -1995
  \multiply\countB by 365
  \advance\countA by \countB
  \countB=\countA
  \divide\countB by 7
  \multiply\countB by 7
  \advance\countA by -\countB
  \advance\countA by 1
  \ifcase\countA\or Sunday\or Monday\or Tuesday\or Wednesday\or
    Thursday\or Friday\or Saturday\fi
\endgroup}

\def\timename{\begingroup
   \countA = \time
   \divide\countA by 60
   \countB = \countA
   \countC = \time
   \multiply\countA by 60
   \advance\countC by -\countA
   \ifnum\countC<10\toks1={0}\else\toks1={}\fi
   \ifnum\countB<12 \toks0={\sevenrm AM}
     \else\toks0={\sevenrm PM}\advance\countB by -12\fi
   \relax\ifnum\countB=0\countB=12\fi
   \hbox{\the\countB:\the\toks1 \the\countC \thinspace \the\toks0}
\endgroup}

\def\timestamp{\dayname, \the\day\ \monthname\ \the\year, \timename}


\print



\def\COMMENT#1\par{\bigskip\hrule\smallskip#1\smallskip\hrule\bigskip}

\def\enma#1{{\ifmmode#1\else$#1$\fi}}

\def\mathbb#1{{\bbold #1}}
\def\mathbf#1{{\bf #1}}


\def\PP{\enma{\mathbb{P}}}


\def\cOO{\enma{\cal O}}


\def\LL{\enma{\mathbf{L}}}

\font\tengoth=eufm10  \font\fivegoth=eufm5
\font\sevengoth=eufm7
\newfam\gothfam  \scriptscriptfont\gothfam=\fivegoth 
\textfont\gothfam=\tengoth \scriptfont\gothfam=\sevengoth
\def\goth{\fam\gothfam\tengoth}
\def \gm {{\goth m}}

\def\GL{\enma{\rm{GL}}}

\def\set#1{\enma{\{#1\}}}


\def\ker{\mathop{\rm ker}\nolimits}

\def\rank{\mathop{\rm rank}\nolimits}

\def\rank{\mathop{\rm rank}\nolimits}
\def\dim{\mathop{\rm dim}\nolimits}

\def\Hom{\mathop{\rm Hom}\nolimits}
\def\Ext{\mathop{\rm Ext}\nolimits}
\def\Tor{\mathop{\rm Tor}\nolimits}

\def\Sym{\mathop{\rm Sym}\nolimits}

\def\Pic{\mathop{\rm Pic}\nolimits}

\def\H{\enma{\rm H}}
%

\newsymbol\boxtimes1202



%
\input diagrams.tex
\diagramstyle[small,labelstyle=\scriptstyle,midshaft]
\overfullrule=0pt

\def\COMMENT#1\par{\bigskip\hrule\smallskip#1\smallskip\hrule\bigskip}

\def\enma#1{{\ifmmode#1\else$#1$\fi}}

\def\mathbb#1{{\bbold #1}}
\def\mathbf#1{{\bf #1}}


\def\PP{\enma{\mathbb{P}}}


\def\cOO{\enma{\cal O}}


\def\LL{\enma{\mathbf{L}}}

\font\tengoth=eufm10  \font\fivegoth=eufm5
\font\sevengoth=eufm7
\newfam\gothfam  \scriptscriptfont\gothfam=\fivegoth 
\textfont\gothfam=\tengoth \scriptfont\gothfam=\sevengoth
\def\goth{\fam\gothfam\tengoth}
\def \gm {{\goth m}}

\def\GL{\enma{\rm{GL}}}

\def\set#1{\enma{\{#1\}}}


\def\ker{\mathop{\rm ker}\nolimits}

\def\rank{\mathop{\rm rank}\nolimits}

\def\rank{\mathop{\rm rank}\nolimits}
\def\dim{\mathop{\rm dim}\nolimits}

\def\Hom{\mathop{\rm Hom}\nolimits}
\def\Ext{\mathop{\rm Ext}\nolimits}
\def\Tor{\mathop{\rm Tor}\nolimits}

\def\Sym{\mathop{\rm Sym}\nolimits}

\def\Pic{\mathop{\rm Pic}\nolimits}

\def\H{\enma{\rm H}}
%

\newsymbol\boxtimes1202

\forward{linear-exactness}{Section}{2}
\forward{complexes}{Section}{3}
\forward{MRC}{Section}{4}
\forward{irred}{Definition}{2.1}
\forward{lexactness}{Proposition}{2.2}
\forward{fail-MRC-curve}{Theorem}{4.1}

\hbox{}
\centerline{\titlefont Exterior algebra methods}
\centerline{\titlefont for the}
\centerline{\titlefont Minimal Resolution Conjecture}
\smallskip
\centerline {by}
\smallskip
\centerline {\bf David Eisenbud, Sorin Popescu, Frank-Olaf
Schreyer, Charles Walter
\footnote{}{\rm The first and the second authors are  grateful
to the NSF for support during the preparation of this work.}
}
\bigskip

\noindent {\bf Abstract:}
{\sl If $r\geq 6, r\neq 9$, we show that the Minimal Resolution
Conjecture fails for a general set of $\gamma$ points in $\PP^r$ for
almost ${1\over 2}\sqrt r$ values of $\gamma$. This strengthens the
result of Eisenbud and Popescu {\rm [1999]}, who found a unique such
$\gamma$ for each $r$ in the given range. Our proof begins like a
variation of that of Eisenbud and Popescu, but uses exterior algebra
methods as explained by Eisenbud and Schreyer {\rm [2000]} to avoid
the degeneration arguments that were the most difficult part of the
Eisenbud-Popescu proof. Analogous techniques show that
the Minimal Resolution Conjecture fails for linearly normal curves of
degree $d$ and genus $g$ when $d\geq 3g-2, \ g\geq4$,
reproving results of 
Schreyer, Green, and Lazarsfeld.}

\section{introd} {Introduction}

\noindent
From the Hilbert function of a
homogeneous ideal $I$ in a polynomial ring $S$ over a field $k$
one can compute a lower bound for the {\it graded betti numbers\/}
$$
\beta_{i,j}= \dim_k\Tor_i^S(S/I,k)_j
$$
because these numbers are the graded ranks of the free modules
in the
minimal free resolution of $I$. The graded betti numbers are
upper semicontinuous in families of ideals with
constant Hilbert function, and in many cases this lower
bound is achieved by the general member of such a family. The
statement that it is achieved for a particular family $T$ is
called the Minimal Resolution Conjecture
(MRC) for $T$, following Lorenzini, who made the conjecture
for the family of ideals of sufficiently general sets of 
$\gamma$ points in $\PP^r$ (all $r,\gamma$).

The minimal resolution conjecture for a general
set of $\gamma$ points in $\PP^r$ has received considerable attention;
see Eisenbud and Popescu [1999] for full references and discussion.
In particular, it is known that the minimal resolution
conjecture is satisfied 
if $r\leq 4$ 
(Gaeta [1951] and [1995], Geramita-Lorenzini [1989], Ballico-Geramita [1986],
Walter [1995], Lauze [1996]), or $\gamma \gg r$
(Hirschowitz-Simpson [1996]), but computer
evidence produced by Schreyer, extended also by 
Mats Boij in his 1994 thesis and by
Beck and Kreuzer in [1996]
suggested that it might fail for certain examples
starting with $11$ points in $\PP^6$.
Indeed, Eisenbud and Popescu
[1999] proved that if $r\geq 6,\ r\neq 9$, then the
minimal resolution conjecture fails for a general set of
$$
\gamma=r+\lfloor(3+\sqrt{8r+1})/2\rfloor
$$
points in $\PP^r$ .
\goodbreak

This family of counterexamples to the minimal resolution conjecture
begins with Schreyer's suggested 11 points in $\PP^6$; and it 
seems to have contained all the examples in print in 1999.
However, the thesis of Boij, completed in 1994 but only
published in [2000], contains computations suggesting
one further counterexample: 21 points in $\PP^{15}$. 

In this paper we will simplify and extend the idea of 
Eisenbud and Popescu to prove the existence
of infinitely many new counterexamples, including the
one suggested by Boij. We prove:

\theorem{failMRC} The Minimal Resolution Conjecture fails for
the general set of $\gamma$ points in $\PP^r$ whenever $r\geq 6,\ 
r\neq 9$, and 
$$
r + 2 + \sqrt{r+2} \; \leq \; \gamma \; \leq \; r + 
(3+\sqrt{8r+1})/2,
$$
and also when $(r,\gamma) = (8,13)$ or $(15,21)$.  

\noindent Thus we get about 
$\left(\sqrt 2 - 1 \right) \sqrt r$ counterexamples for each $r\geq
6,\ r\neq 9$.

\smallskip
Our proof not only gives more, but it avoids the subtle degeneration
argument used by Eisenbud and Popescu. In the presentation
below we will skip some of the steps presented in 
Eisenbud-Popescu [1999], and fully treat only the 
new ideas, so it may be helpful to 
the reader to explain their strategy and the point at which
it differs from
ours. Their
proof can be divided into three steps, starting from a set
$\Gamma$ of $\gamma$ general points in $\PP^r$:

\medskip
\item{I.} 
Their crucial first step is to consider the Gale transform of
$\Gamma\subset\PP^r$, which is a set $\Gamma'$ of $\gamma$ points in
$\PP^s$ with $\gamma=r+s+2$ determined ``naively'' as follows: If the
columns of the $(r+1) \times \gamma$ matrix $M$ represent the points
of $\Gamma$, then a $(s+1) \times \gamma$ matrix $N$ representing
$\Gamma'$ is given by the transpose of the kernel of $M$ (see
Eisenbud-Popescu [1999] for a precise definition).  The free
resolution of the canonical modules of the homogeneous coordinate
rings of $\Gamma$ and $\Gamma'$ are related, and using this relation
and the fact that (under our numerical hypotheses) the ideal of
$\Gamma'$ is not $2$-regular, they construct a map $\phi_\Gamma$ from
a certain linear free complex $F_\bullet(\mu)$ to the dual of the
resolution of the ideal $I_\Gamma$ of $\Gamma$.  For the $r$ and $\gamma$
considered in the theorem, straightforward arithmetic shows that the
injectivity of the top degree component of $\phi_\Gamma$ would force
the graded betti numbers of $I_\Gamma$ to be too large for the minimal
resolution conjecture to hold. (We will recall the definition of
$F_\bullet(\mu)$ in \ref{complexes} below.)
\smallskip

\item{II.} They show that $\phi_\Gamma$ is an injection of complexes
if $F_\bullet(\mu)$ has a property somewhat weaker than
exactness called {\it linear exactness\/} or 
{\it irredundancy} (see \ref{irred} below).
\smallskip

\item{III.}
For a generic set of points $\Gamma$ they show that $F_\bullet(\mu)$
is irredundant by a subtle degeneration argument:
they degenerate the general set $\Gamma$ to
lie on a special curve $C$, and the argument finishes with a study of
a refined stability property of the tangent bundle of the projective
space restricted to $C$ in a certain embedding of $C$
connected with the Gale transform of the points.
\medskip

In our new proof step I remains unchanged, but II is
replaced by a more precise statement, which requires us
to verify a weaker condition in III; this weak condition
can be verified without any degeneration argument, making
the replacement for step III much simpler. More precisely
the injectivity of
$\phi_\Gamma$ is only needed for a sufficiently large irredundant quotient of
the complex $F_\bullet(\mu)$, and we give a criterion for this
weaker condition using exterior algebra methods.  The new criterion is
checked in step III with
an argument inspired by Mark Green's proof [1999] of the Linear Syzygy
Conjecture.  The generality of $\Gamma$ enters via work of
Kreuzer [1994] showing that certain multiplication maps of the
canonical module of the cone over $\Gamma$ are 1-generic in the sense
of Eisenbud [1988].

In \ref{linear-exactness} we introduce and study the irredundancy of
linear complexes using the Bernstein-Gel'fand-Gel'fand (BGG)
correspondence. In \ref{complexes} we construct the complex whose
irredundancy is the key to the failure of the minimal resolution
conjecture. In the final \ref{MRC} we briefly explain how Gale duality
allows this theory to be applied to sets of points, and give the
arithmetic part of the proof. As another application of our 
techniques we recover the result on curves
of Schreyer (unpublished) and Green-Lazarsfeld [1988]; see \ref{fail-MRC-curve}.

It is a pleasure to thank Tony Iarrobino and the Mathematisches
Forschungsinstitut Oberwolfach.  Iarrobino, while reviewing Boij
[2000] for MathSci, noticed that Boij's example of $21$ points in
$\PP^{15}$ was not in the Eisenbud-Popescu list.  The authors of the
present paper were at that moment attending a conference at Oberwolfach.
Iarrobino triggered our collaboration by asking Eisenbud and Popescu
whether their methods could encompass the new example.  The joy and
success we had in answering his question owe much to the marvelous
facilities and atmosphere in the Black Forest!
\medskip

\noindent{\bf Notation:} Throughout this paper $k$ denotes an arbitrary field.
Let $V$ be a finite-dimensional vector space over $k$.
We work with graded modules over $S:=\Sym(V)$ and $E:=\wedge V^*$.
We think of elements of $V$ as having degree 1 and
elements of $V^*$ as having degree $-1$.
We write $\gm$ for the maximal ideal generated by $V$ in $S$.

\section{linear-exactness} {Linear complexes and exterior modules}

In this section we illustrate the exterior algebra approach to
linear free resolutions with some basic ideas used
in the rest of the paper, and with a new proof of
the linear rigidity theorem of Eisenbud and Popescu [1999].

\definition{irred} A complex of graded free $S$-modules 
$$
{F_\bullet}:
\quad \cdots
\rTo^{\phi_{3}} F_{2}
\rTo^{\phi_{2}} F_{1} 
\rTo^{\phi_{1}} F_{0}
$$
is a called a {\rm standard linear free complex\/} if $F_i$ is generated
in degree $i$ for all $i$ (note that all the differentials $\phi_i$
are then represented by matrices of linear forms.) \hfill\break 
A {\rm linear free complex\/} is any twisted (or shifted) standard linear
free complex (that is, a right bounded complex of free modules $F_i$
so that for some integer $s$ the generators of $F_i$ are all in degree
$s+i$ for all $i$).
\hfill\break
We say that a linear free complex
$F_\bullet$ is {\rm irredundant\/} (Eisenbud-Schreyer {\rm [2000]}) or 
{\rm linearly exact} (Eisenbud-Popescu {\rm [1999]}) if
the induced maps $$F_{i+1}/\gm F_{i+1}\to \gm F_i/\gm^2 F_i$$
are monomorphisms for $i>0$.

One sees immediately that a standard linear free complex $F_\bullet$
as above is irredundant if and only if $\H_i(F_\bullet)_i=0$ for all
$i > 0$.  

As explained in Eisenbud-Schreyer [2000],
irredundancy can be conveniently rephrased via the
Gel'fand-Gel'fand-Bernstein correspondence as follows:

If $P$ is a graded $E$-module we define a linear free $S$-complex
$F_\bullet=\LL(P)$ with free modules $F_i=S\otimes_kP_i$ 
generated in degree $i$ and differentials
$$
\phi_i:\ F_i\to F_{i-1}\qquad \qquad
1\otimes p\mapsto \sum x_i\otimes e_ip\in S\otimes P_{i-1},
$$
where $\{x_i\}$ and $\{e_i\}$ are (fixed) dual bases of $V$ and $V^*$.
Up to twisting and shifting every linear free complex of
$S$-modules arises from a unique $E$-module in this way. The
linear complex $\LL(P)$ is standard if $P_i=0$ for all $i<0$.

Given a graded $E$-module $P$ we write $P^*$ for the dual graded
$E$-module $P^*=\Hom_k(P,k)$. We take $k$ in degree $0$,
so that $(P^*)_i$ is the dual vector space to $P_{-i}$.

\proposition{lexactness} If 
$$
F_\bullet=\LL(P): \quad \cdots \rTo F_2\rTo F_1\rTo F_0
$$
is a standard linear free complex of $S$-modules, then $F_\bullet$ is
irredundant if and only if $P^*$ is generated
as an $E$-module in degree $0$. 
\Box

Thus any linear free complex $F_\bullet=\LL(P)$ as in \ref{lexactness}
has a unique maximal irredundant quotient $F'_\bullet$, which is
functorial in $F_\bullet$, constructed as follows: If $Q$ is the
$E$-submodule of $P^*$ generated by $P_0^*$, then
$F'_\bullet:=\LL(Q^*)$.  Alternatively, $F'_\bullet=\LL(P/N)$, where
$N$ is the submodule with graded pieces $N_i = \set{n \in P_i \mid
\wedge^i V^* \cdot n = 0}$.  Observe that we have $F'_0=F_0$.

The next lemma is the basic tool which will allow us to deduce that a
minimal free resolution is large.  A weaker version of the result was
implicit in Eisenbud-Popescu [1999].

\lemma{inclusion} Suppose that 
$
\alpha_\bullet:{F_\bullet}\to G_\bullet
$
is a map from an irredundant standard linear free complex to a
minimal free complex.
If $\alpha_0$ is a split monomorphism, then $\alpha_i$ is a
split monomorphism for all $i\geq 0$.

\proof By induction it suffices to prove that $\alpha_1$ is
a split monomorphism.  Since both $F_\bullet$ and
$G_\bullet$ are minimal (that is, the differentials are
represented by matrices of elements of $\gm$) there is a commutative
diagram
$$\diagram
F_1/\gm F_1&\rTo&\gm F_0/\gm^2 F_0\cr
\dTo^{\overline\alpha_1}&&\dTo\cr
G_1/\gm G_1&\rTo&\gm G_0/\gm^2 G_0
\enddiagram
$$
whose vertical maps are induced by $\alpha_\bullet$ and whose 
horizontal maps are induced by the differentials of $F_\bullet$ and
$G_\bullet$. 

Since $\alpha_0$ makes $F_0$ a summand of $G_0$, the right hand vertical map
is  a monomorphism. Because $F_\bullet$
is irredundant, the top map is a monomorphism. It follows that the map 
$\overline\alpha_1:F_1/\gm F_1 \to G_1/\gm G_1$ induced by $\alpha_1$
is a monomorphism, and since $G_1$ is free, this implies that $\alpha_1$
is split.\Box

\ref{inclusion} may be applied to give a lower bound on
betti numbers:

\proposition{superfunctoriality} Suppose that 
$\alpha_\bullet: F_\bullet\to G_\bullet$ is a map from a standard
linear complex to a minimal free resolution. There exists a map
$\beta_\bullet: F_\bullet\to G_\bullet$ which is homotopic to
$\alpha_\bullet$ and such that $\beta_\bullet$ factors through the
maximal irredundant quotient $F'_\bullet$ of $F_\bullet$.  Further, if
$\alpha_0: F_0\to G_0$ is a split monomorphism, then the rank of $G_i$ is at
least as big as the rank of $F'_i$, for all $i$.

\proof Let $\gamma_\bullet:F_\bullet\to F'_\bullet$ be the natural
map to the maximal irredundant quotient. The image of $F'_1$ in $F'_0=F_0$ is
contained in that of $F_1$, so there is a map $\beta_1: F'_1\to G_1$
lifting
$\alpha_0: F'_0\to G_0$. Since $G_\bullet$ is acyclic we may continue to
lift, and so inductively we get a map of
complexes $\beta'_\bullet: F'_\bullet \to G_\bullet$. We take
$\beta_\bullet = \beta'_\bullet\gamma_\bullet$ and thus we have
$\beta_0=\alpha_0$. Since $F_\bullet$ is a free complex and
$G_\bullet$ is acyclic, this implies that $\alpha_\bullet$ is
homotopic to $\beta_\bullet$.

The second statement follows by applying \ref{inclusion}
to the map $\beta'_\bullet$.
\Box

In Eisenbud-Popescu [1999] a rigidity result for irredundancy is
required. The proof given there is just a reference to the result on
the rigidity of Tor proved by Auslander and Buchsbaum in [1958]. 
The exterior method yields a novel proof of this result (which
we will not need in the sequel):

\proposition{general-rigidity} Let $R$ be a graded commutative
or anti-commutative ring, and let $M$ be a graded $R$-module
which is generated by $M_0$. If $R'\subset R$ is a graded
subring such that $R'_1\cdot M_0=M_1$, then $M$ is generated
by $M_0$ also as an $R'$-module.

\proof We show by induction on $n$ that $R'_1\cdot M_n=M_{n+1}$,
the initial case $n=0$ being the hypothesis. If $R'_1\cdot M_{n-1}=
M_{n}$, then $R'_1\cdot M_n=R'_1\cdot (R_1M_{n-1})=
R_1\cdot (R'_1M_{n-1})=R_1\cdot M_n=M_{n+1}$, where the
second equality holds  by the (anti)-commutativity assumption
and the third by the induction hypothesis.\Box

\corollary{linear-rigidity} (Linear Rigidity) 
Let $S=k[x_0,\dots,x_r]$
be a polynomial ring, and let 
$$
F_\bullet: \quad \cdots \rTo F_2
\rTo F_1\rTo F_0
$$ 
be an irredundant standard linear free complex.
Let $S'=S/I$ be a graded quotient of $S$.
The complex $G_\bullet:=S'\otimes_S F_\bullet$ is irredundant
iff $\H_{1}(G_\bullet)$ is zero in degree $1$.
 
\proof The irredundancy of $G_\bullet$ involves only information about 
$S'_0$ and $S'_1$, so we may assume that $S'$ is a polynomial ring
$S'=\Sym(V')$ with $V'=V/I_1$. Write $E'=\wedge(V'^*)\subset
E=\wedge(V^*)$ for the corresponding exterior algebras. If
$F_\bullet=\LL(P)$, then the complex $G_\bullet$ over $S'$ corresponds
to the same graded vector space $P$, regarded as an $E'$-module by
restriction of scalars.  Consider now the truncated complex
$$
K_\bullet : \quad 
\cdots \rTo 0 \rTo 0 \rTo G_1 \rTo G_0.
$$
The hypothesis $\H_1(G_\bullet)_1 = 0$ implies that $\H_1(K_\bullet)_1
= 0$ and hence, by \ref{lexactness}, that $P_1^* =V'^*\cdot P_0^*$. On
the other hand, $P^*$ is generated as an $E$-module by $P_0^*$ since
$F_\bullet$ is an irredundant standard linear complex.  So by
\ref{general-rigidity} we have also $P^* =E' \cdot P_0^*$.  By
\ref{lexactness} this suffices.\Box

\remark{} A similar ``linear rigidity'' result for complexes
over the exterior algebra can be also deduced from
\ref{general-rigidity} applied this time for $R$ a
polynomial ring.

\section{complexes} The complexes $F_\bullet(\mu)$

\medskip
We will apply \ref{lexactness} and \ref{inclusion} to a family of
linear complexes $F_\bullet(\mu)$ (which were called $E^{-1}(\mu)$ in
Eisenbud-Popescu [1999]).  It is convenient here to define them in
terms of the functor $\LL$.

Let $U$, $V$, and $W$ be finite-dimensional vector spaces of
dimensions $u$, $v$, and $w$, respectively.  Let $A := \sum_l
\wedge^{l} W^*\otimes \Sym_l(U^*),$ and let $Q := \sum_l
\wedge^{l+1}W^*\otimes \Sym_l(U^*)$.  Then $A$ is an anticommutative
graded algebra, and $Q$ is a graded $A$-module.  If $\mu : W \otimes U
\to V$ is a pairing, then the dual $\mu^*:V^*\to W^*\otimes U^*$
extends to a map of algebras $\tilde \mu: E=\wedge(V^*)\to A$.  We
regard $Q$ as an $E$-module via $\tilde\mu$, and we consider the corresponding
linear free complex $F_\bullet(\mu) := \LL(Q^*)$ over $S=\Sym(V)$:
$$
F_\bullet(\mu): \quad 
0 \to \wedge^w W\otimes D_{w-1}(U) \otimes S(-w+1) \to \cdots \to 
\wedge^2 W \otimes U \otimes S(-1) \to W \otimes S,
$$
where $D_m(U)$ denotes the $m$-graded piece of the divided power algebra.

Though we shall not need the formula, it is easy to 
give the differentials explicitly:
$
\delta_{l}(\mu): F_{l}(\mu) \to F_{l-1}(\mu),
$
is the composite of the tensor product of the diagonal maps
of the exterior and divided powers algebras
$$
\wedge^{l+1}W\otimes D_{l}U\otimes S(-l)\rTo
\wedge^{l}W\otimes W\otimes D_{l-1}U\otimes U\otimes S(-l),
$$
and the map induced by the pairing $\mu$
$$
\wedge^{l}W\otimes W\otimes D_{l-1}U\otimes U\otimes S(-l)
\rTo
\wedge^{l}W\otimes D_{l-1}U\otimes S(-l+1).
$$

Note that if $V=W\otimes U$ and $\mu$ is the identity map, then $Q$ is
generated as an $E$-module by $Q_0$, so the complex $F_\bullet(\mu)$
is irredundant by \ref{lexactness}.  The proof in Eisenbud-Popescu
[1999] shows that $F_\bullet(\mu)$ remains exact when $\mu$ is
specialized to a certain pairing coming from the canonical module of a
generic set of $\gamma$ points in $\PP^r$ for suitable $\gamma,
r$. Here we are instead interested in knowing whether $F_\bullet(\mu)$
has an irredundant quotient complex with the same first and last terms
as $F_\bullet(\mu)$. This condition turns out to be a lot easier to analyze!

Following Eisenbud [1988] we say that $\mu : W \otimes U \to V$ is
{\it $1$-generic\/} if $\mu(a\otimes b)\neq 0$ for all nonzero vectors
$a \in W$ and $b \in U$. We say that $\mu$ is {\it geometrically
$1$-generic\/} if the induced pairing $\mu\otimes\overline k$ is
$1$-generic, where $\overline k$ is the algebraic closure of $k$.  If
$\mu$ is geometrically $1$-generic, then by a classic and elementary
result of Hopf the dimensions $u,v,w$ of the three vector spaces
satisfy the inequality $v\geq u+w-1$.

The main result of this section is the following:

\theorem{1-generic}  If the pairing $\mu:\; W\otimes U\to V$ is
geometrically $1$-generic, and $W \neq 0$ and $U \neq 0$, then
$F_\bullet(\mu)$ and its maximal irredundant quotient
$F'_\bullet(\mu)$ have the same last term $F_{w-1} = F'_{w-1}$,
that is 
$$
{F'_\bullet(\mu)}:\quad \wedge^w W\otimes D_{w-1}(U)\otimes S(-w+1)
\rTo F'_{w-2} \rTo \cdots 
\rTo F'_{1} \rTo W\otimes S.
$$

\proof We use notation as above. Let $Q'$ be the submodule of $Q$ generated by
$Q_0$, and set $P=Q^*, P'=Q'^*$, 
so that $F_\bullet=\LL(P)$ and
$F'_\bullet=\LL(P')$. 
By \ref{lexactness}, the conclusion of the
Theorem is equivalent to the statement that the $E$-multiplication map
\newarrow{Equals}=====
$$
\diagram[small,midshaft]
m: & E_{w-1}\otimes_k P_0 &\rTo & P_{-w+1}\\
& \dEquals & &  \dEquals\\
&\wedge^{w-1} V^*\otimes W^*& \rTo& \wedge^w W^*\otimes \Sym_{w-1}(U^*)\\
\enddiagram
$$
is surjective. Via the natural identification
$W^*\otimes \wedge^w W\cong\wedge^{w-1}W$
this is equivalent to the map
$$
m': \wedge^{w-1} V^*\otimes_k \wedge^{w-1}W\rTo\Sym_{w-1}(U^*)
$$
induced by $\mu$ being surjective. A straightforward computation shows
that the image of $m'$ is the space generated
by the $(w-1)\times (w-1)$ minors of the linear map of free
$\Sym(U^*)$-modules
$$
\bar\mu: \quad W\otimes\Sym(U^*)\rTo V\otimes\Sym(U^*)(1)
$$
associated to $\mu$, and the next lemma shows that under the
$1$-genericity assumption these minors span $\Sym_{w-1} (U^*)$.  
This proves the theorem.
\Box

\lemma{minors}
If $\mu : W \otimes U \to V$ is a geometrically $1$-generic pairing of
finite-dimensional vector spaces, and if $d$ is an integer such that $d \leq \dim W$, then
the $d \times d$ minors of the associated linear map $\bar \mu : W \otimes
\Sym(U^*) \to V \otimes \Sym(U^*)(1)$ of free
$\Sym(U^*)$-modules span $\Sym_{d} (U^*)$.

\proof  
It is enough to prove the lemma for an algebraically closed base field
$k$, so that the Nullstellensatz holds.  Also, $u,v,w$ continue to
be the dimensions of the three vector spaces.

We first show that, for an arbitrary subspace $W'\subset W$ of
dimension $d$, the $d\times d$ minors of the restricted
map
$$
\bar\mu': W'\otimes\Sym(U^*)\rTo V\otimes\Sym(U^*)(1),
$$
regarded as polynomial functions on $U$, have a common zero only at
the origin $0 \in U$.  For $\bar\mu'$ is essentially a $v \times d$
matrix whose entries are linear functions on $U$, so its $d
\times d$ minors have a common zero at a point $b \in U$ if and only
if the matrix's columns, when evaluated at $b$, are not linearly
independent.  But this means that there is a $0
\neq a \in W'$ such that $\mu(a
\otimes b) = 0$, which gives $b = 0$ by $1$-genericity.

A general position argument shows that, for a general projection
$V\rOnto V'$ with $\dim V'= d+u-1$, the only common zero of the
$d\times d$ minors of the composite map
$$
\bar\mu'': W'\otimes\Sym(U^*)\rTo V'\otimes\Sym(U^*)(1),
$$
is still $0 \in U$.  
(One may argue that the dual of $\bar\mu'$ may be regarded as
a surjective map $\cOO^v\rOnto \cOO^{d}(1)$ 
of vector bundles on $\PP^{u-1}$, and a rank $d$ vector
bundle generated by global sections on a $(u-1)$ dimensional space
can always be generated by just $d+u-1$ of them; 
see for example Eisenbud-Evans [1973], or Serre [1958]).

The expected codimension of the locus defined by the $d \times d$
minors of the $(d+u-1) \times d$ matrix representing $\bar\mu''$ is $u$, which is
equal to the actual codimension. Thus the Eagon-Northcott complex
resolving the minors of $\bar\mu''$ is exact, and it follows that
there are  ${d+u-1\choose d}$ linearly independent such minors.
Since this is also the dimension of $\Sym_{d}(U^*)$, we see that the
minors span $\Sym_{d}(U^*)$ as required.\Box

We say that a map of complexes 
$\alpha_\bullet:\ F_\bullet\to G_\bullet$
is a {\it degreewise split injection\/} if each $\alpha_i: F_i\to G_i$
is a split injection (this is weaker than being a split injection
of complexes); in this case we say that $F_\bullet$ is a
{\it degreewise direct summand\/} of $G_\bullet$.
\smallskip

The irredundant complexes
$F'_\bullet(\mu)$ arise
in geometric situations as follows. 

\theorem{MRC-fail} Let $L,L'$ be two line bundles
on a scheme $X$ over a field $k$, and let 
$L'' := L \otimes L'$.  Let
$W\subset \H^0(L),\ U \subset \H^0(L')$
be nonzero finite-dimensional linear series such that the multiplication 
$$\mu: W\otimes U\rTo \H^0(L'')$$ is geometrically $1$-generic.
Let
$V \subset \H^0(L'')$ be a finite-dimensional linear series
containing $W\cdot U$. 
Let $S:=\Sym(V)$, and let $M \subset \oplus_{n\geq 0} \H^0(L \otimes
{L''}^{\otimes n})$ be a finitely generated graded $S$-submodule such
that $W \subset M_0$.  Then the
maximal irredundant quotient 
$$
{F'_\bullet(\mu)}:\quad \wedge^w W\otimes D_{w-1}(U)\otimes S(-w+1)
\rTo F'_{w-2} \rTo \cdots 
\rTo F'_{1} \rTo W\otimes S
$$
of the linear free complex $F_\bullet(\mu)$ injects as a degreewise
direct summand of the minimal free resolution of $M$. (Here as above
$w=\dim W$.)

\proof
Let $G_\bullet$ be the minimal free resolution of $M$ over $S$.  The
hypotheses on $M$ imply that $G_0 = (W \otimes S) \oplus L$, with $L$
a free $S$-module.  We can construct inductively a morphism
$\alpha_\bullet : F_\bullet(\mu) \to G_\bullet$, 
such that $\overline{\alpha_0}:F_0(\mu)/\gm F_0(\mu)
\to G_0/\gm G_0$ is a monomorphism, as a lifting
\diagram[small]
\cdots & \rTo & \wedge^3 W \otimes D_2 U \otimes S(-2) & \rTo &
\wedge^2 W \otimes U \otimes S(-1) & \rTo & W \otimes S \\
&& \dDotsto && \dDotsto^{\alpha_1} && \dInto \\
\cdots & \rTo & G_2 & \rTo & G_1 & \rTo & (W \otimes S) \oplus L &
\rOnto & M 
\enddiagram
For instance, $\alpha_1$ exists since 
the composition $\wedge^2 W \otimes U \otimes S(-1) \to W
\otimes S \to M$ vanishes.
So by \ref{inclusion}, $\alpha_\bullet$ induces a degreewise split
inclusion $F'_\bullet(\mu) \rInto G_\bullet$, with $F'_\bullet(\mu)$
having the form given in \ref{1-generic} (since $\mu$ is geometrically
$1$-generic).
\Box

It is easy to see that if $X$ is a geometrically reduced and irreducible
scheme, then any pairing $\mu: W\otimes U\rTo V$ induced
by multiplication of sections as in \ref{MRC-fail}
is geometrically $1$-generic. The hypothesis sometimes holds
for reducible or nonreduced schemes too, as in the
case we will use for the proof of \ref{failMRC}.

The next example shows that the complexes $F_\bullet(\mu)$ and
$F'_\bullet(\mu)$ are not always the same.

\example{cex-eb-linearly-exact} 
If we take $U=W=\H^0(\cOO_{\PP^1}(2))$ with
$\mu$ the multiplication map to $V=\H^0(\cOO_{\PP^1}(4))$, that is
$$
\mu: {k[s,t]}_2\otimes {k[s,t]}_2=W\otimes U\rTo V={k[s,t]}_4,
$$
then $\mu$ is
geometrically $1$-generic because $k[s,t]$ is a domain
and remains so over the algebraic closure of $k$.
The module $M=\oplus_{n\ge 0}\H^0(\cOO_{\PP^1}(2)\otimes\cOO_{\PP^1}(4n))$ 
of \ref{MRC-fail} is the graded module
associated to the line bundle 
$\cOO_{\PP^1}(2\hbox{ points})$ on $\PP^1$,
regarded as a module over the homogeneous coordinate ring
$S$ of $\PP^4$ via the embedding of $\PP^1$ as the rational
normal quartic $C$. This line bundle is $\omega_C(1)$, the twist of the
canonical line bundle by the hyperplane line bundle, so 
the minimal resolution of $M$ is the dual
of the minimal resolution of the homogeneous coordinate
ring of $C$, suitably shifted. It has the form
$$
G_\bullet: \quad 0\rTo S(-4)\rTo S^6(-2)\rTo S^8(-1)\rTo S^3\rTo M\rTo 0.
$$
On the other hand, $F_\bullet$ has the form
$$F_\bullet(\mu):\quad
0\rTo S^6(-2)\rTo (S^3\otimes S^3)(-1)\rTo^{\delta_1(\mu)} S^3,
$$
and $F_1=S^9(-1)$ cannot inject into $G_1=S^8(-1)$. The columns
of the matrix of $\delta_1(\mu)$
have exactly one $k$-linear dependence relation, so 
the maximal irredundant quotient of $F_\bullet$ has the form
$$
F'_\bullet: \quad
0\rTo S^6(-2)\rTo S^8(-1)\rTo S^3,
$$
so in particular $F'_2=F_2$, in accordance with \ref{1-generic}.
By \ref{inclusion} any map from $F_\bullet$ to $G_\bullet$
lifting the identity on $F_0$ is a degreewise split monomorphism
on $F'_\bullet$; in fact $F'_\bullet$ 
is isomorphic to the linear strand of $G_\bullet$ in this case.

\section{MRC} Failure of the Minimal Resolution Conjecture

Before proving \ref{failMRC}, we present the bare bones of Gale
duality.  The reader who wishes to go further into this rich and
beautiful subject can consult Eisenbud-Popescu [1999, 2000],
or Dolgachev-Ortland [1988].

Let $\Gamma \subset \PP^r$ be a set of $\gamma = r+s+2$ points.  The
linear forms on $\PP^r$ define a linear series $V \subset
\H^0(\cOO_\Gamma(1))$.  The orthogonal complement $V^\perp \subset
\H^0(\cOO_\Gamma(1))^*$ is identified by Serre duality with a
linear series $(\omega_\Gamma)_{-1} \subset \H^0(K_\Gamma(-1))$
(where $K_\Gamma$ denotes the canonical sheaf of $\Gamma$). 
Under mild hypotheses on $\Gamma$, the linear series
$(\omega_\Gamma)_{-1}$ is very ample and gives an embedding whose
image $\Gamma' \subset \PP^s$ is the {\it Gale transform} of $\Gamma$.
As a set it is well-defined only 
modulo linear changes of coordinates on $\PP^s$ and $\PP^r$.

Concretely, if the columns of the $(r+1) \times \gamma$ matrix $M$
represent the points of $\Gamma\subset\PP^r$, and if $N$ is a $\gamma \times
(s+1)$ matrix whose columns are a basis of $\ker(M)$, then the rows of
$N$ represent the points of $\Gamma'\subset\PP^s$.

In the proof of \ref{failMRC} we need four facts about the
Gale transform, which the reader may easily verify (or find
in Eisenbud-Popescu [1999, 2000]).

\medskip
\item{--}
{\sl The subset $\Gamma' \subset \PP^s$ is defined for a
general $\Gamma$}.  In fact, it is defined as long as no subset containing all
but two of the points of $\Gamma$ lies in a hyperplane of $\PP^r$.

\item{--}
{\sl The Gale transform of $\Gamma'$ is $\Gamma$}.

\item{--} 
{\sl There is a natural identification of $V = \H^0(\cOO_{\PP^r}(1))$
with $(\omega_{\Gamma'})_{-1}$}. 

\item{--}
{\sl $\Gamma$ is general if and only if $\Gamma'$ is general}.  More
formally, a $\GL(s+1)$-invariant
Zariski open condition on $\Gamma'$ induces a $\GL(r+1)$-invariant
Zariski open condition on $\Gamma$.

\medskip
\noindent{\sl Proof of \ref{failMRC}.\/}
Let $s\geq 3$ be an integer and set
$r={s+1\choose 2}+\delta$. Suppose that
$$
0\le \delta\le {s\choose 2}-\cases{1 &if $s\leq 4$,\cr 2&if $s \geq 5$,}
$$
and let
$\Gamma\subset\PP^r$ be a set of $\gamma(r,s)=r+s+2$ points,
in linearly general position. We shall show that the
Minimal Resolution Conjecture fails for $\Gamma$.
From this statement it is
easy to solve for $\gamma$ in terms of $r$, obtaining the
range given in \ref{failMRC}.

Let $I_\Gamma$ be the
homogeneous ideal of $\Gamma\subset\PP^r$ and let
$$
\omega_{\Gamma} = \Ext^{r-1}_S(I_\Gamma, S(-r-1))
$$
be the canonical module. The free resolution of $\omega_\Gamma$ is, up
to a shift in degree, the dual of the resolution of $S/I_\Gamma$.

We will assume that $\Gamma$ is a general set of points such
that: its Gale transform $\Gamma' \subset \PP^s$ is defined and
is in linearly general position, $\Gamma$ imposes independent
conditions on quadrics, and $\Gamma'$ is not contained in any
quadric (note that ${s+2\choose 2} <
\gamma(r,s) \leq {r+2\choose 2}$).

Since $\gamma(r,s)>{s+2\choose 2}$, it follows that
$\Gamma'$ does not impose independent conditions on quadrics
of $\PP^s$, and thus $\H^1({\cal I}_{\Gamma'}(2))\neq 0$.
We set
$$
U := (\omega_{\Gamma'})_{-2}=
\Hom_k(\H^1({\cal I}_{\Gamma'}(2)), k)
\neq 0.
$$
(In this setting we have $\dim U = \delta+1$.)  Write $W = \H^0({\cal
O}_{\PP^s}(1))$, and $V=\H^0({\cal O}_{\PP^r}(1))$.  Multiplication
induces a natural pairing
$$
\mu:\quad W\otimes U \rTo V=(\omega_{\Gamma'})_{-1}.
$$

We now carry out the checks needed to apply \ref{MRC-fail}.  First,
$\Gamma'\subset\PP^s$ is in linearly general position, so by Kreuzer
[1994] the pairing $\mu$ is $1$-generic over any field and thus
geometrically $1$-generic.  Second, $\Gamma \subset \mathbb P^r$
imposes independent conditions on quadrics, so $(\omega_\Gamma)_{-d} =
0$ for $d \geq 2$.  Hence $\omega_\Gamma(-1)$ is a finitely generated
graded $S$-submodule of $\oplus_{n \geq 0} \H^0(K_\Gamma(n-1))$
containing $W = (\omega_\Gamma)_{-1}$ in degree $0$.

Therefore the irredundant linear free quotient $F'_\bullet(\mu)$ of
the complex $F_\bullet(\mu)$ injects onto a degreewise split direct
summand of the minimal free resolution of $\omega_\Gamma(-1)$, which
is the dual of the minimal free resolution of $(S/I_\Gamma)(r+2)$.
Hence the true graded betti number $\beta_{(r-s),(r-s+2)}$ for
$S/I_\Gamma$ satisfies
$$
\beta_{(r-s),(r-s+2)} \geq \rank F'_{s}(\mu) = {s+\delta \choose \delta}>0.
$$ 
The corresponding expected graded betti number
$\widetilde\beta_{(r-s),(r-s+2)}$ is found in Eisenbud-Popescu [1999],
or it can be computed by taking the Hilbert series for $\Gamma$
$$
\sum_i \dim (S/I_\Gamma)_i \cdot t^i = 1 + (r+1)t + 
\gamma(r,s) \cdot {{t^2}\over{1-t}} = 
{{\sum_j b_j t^j}\over{(1-t)^{r+1}}}
$$
and calculating that 
$$
(-1)^{r-s} b_{r-s+2} = 
{{(2\delta+4-s^2+s)}\over{(s^2-s+2\delta+4)}}\cdot 
{{s+1\choose 2}+\delta\choose s}.
$$
This number is the alternating sum of the graded betti numbers of
degree $r-s+2$, in this case $\beta_{(r-s),(r-s+2)} -
\beta_{(r-s+1),(r-s+2)}$.  The expected graded betti number is thus
$$
\widetilde\beta_{(r-s),(r-s+2)} = \max \left( 
(-1)^{r-s} b_{r-s+2}, 0 \right).
$$
The MRC certainly fails when $\widetilde\beta_{(r-s),(r-s+2)} = 0$,
which is equivalent to $0 \leq \delta \leq {s \choose 2} - 2$.  
The MRC also fails for $s = 3,4$ and $\delta = {s \choose 2} - 1$
because ${ s + \delta \choose \delta} >
\widetilde\beta_{(r-s),(r-s+2)}$ in those cases.  (For $s=4$ this is
Boij's example!)  Straightforward computation shows that there are no
other cases where ${ s + \delta \choose \delta} >
\widetilde\beta_{(r-s),(r-s+2)}$. 
\Box

Schreyer observed around 1983 (unpublished) that the Minimal Resolution
Conjecture fails for linearly normal curves of degree $d>g^2-g, \ g\geq 4$.
A result of Green and Lazarsfeld [1988] implies that this bound
can be improved to $d>3g-3$. We can use the technique developed
above to give a new proof of the relevant part of the Green-Lazarsfeld
result, recovering the failure of the MRC in a new way.

\theorem{fail-MRC-curve} Let $C$ be a smooth curve of genus $g\ge 2$,
let $L\in\Pic^d(C)$ be a line bundle on $C$ of degree $d\ge 2g+2$, and
denote by $\varphi_{|L|}:C\to \PP^{d-g}$ the embedding defined by the
complete linear system $|L|$.
\item{$a)$} If $\H^0(L\otimes\omega_C^{-1})\neq 0$ (that is, equivalently, 
if the curve $\varphi_{|L|}(C)\subset \PP^{d-g}$ has a $(d-2g+1)$-secant
$\PP^{d-2g-1}$), then the maximal irredundant quotient of the
free linear complex
$$
{F_\bullet}:=\LL(\oplus_{l\ge 0}\wedge^{l+1}\H^0(\omega_C)\otimes
D_l(\H^0(L\otimes\omega_C^{-1}))
$$ 
injects as a degreewise direct summand of the minimal free resolution
of the $\Sym(\H^0(L))$-module $\oplus_{n\ge 0}\H^0(\omega_C\otimes
L^n)$. In particular $\varphi_{|L|}(C)\subset \PP^{d-g}$
doesn't satisfy the property $N_{d-2g}$. 
\item{$b)$} If $g\ge 4$ and $d> 3g-3$,  then 
the Minimal Resolution Conjecture fails for $\varphi_{|L|}(C)\subset
\PP^{d-g}$.

\remark{Np} By Green [1984],  $\varphi_{|L|}(C)\subset \PP^{d-g}$ satisfies the
property $N_{d-2g-1}$ (that is the embedding is projectively normal, the homogeneous ideal 
$I_{\varphi_{|L|}(C)}$ is generated by quadrics and all syzygies in
the first $(d-2g-2)$ steps in its minimal free resolution are
linear). On the other hand, if $C$ is non-hyperelliptic then property
$N_{d-2g}$ fails for $\varphi_{|L|}(C)\subset \PP^{d-g}$ iff the curve
has a $(d-2g+1)$-secant $(d-2g-1)$-plane, by Green-Lazarsfeld [1988,
Theorem 2]. \ref{fail-MRC-curve}, $a)$  spells out explicitly this
failure!
\medskip

\noindent{\it Proof of \ref{fail-MRC-curve}}: 
By hypothesis $d\ge 2g+2$ so 
$\varphi_{|L|}(C)\subset \PP^{d-g}$ is arithmetically
Cohen-Macaulay, its homogeneous ideal $I_{\varphi_{|L|}(C)}$ 
is generated by quadrics and is $2$-regular. 
Now \ref{MRC-fail}, applied to the geometrically $1$-generic pairing
$$
\mu: \H^0(\omega_C)\otimes\H^0(L\otimes\omega_C^{-1})\to\H^0(L)
$$  
yields the complex $F_\bullet=F_\bullet(\mu)$ whose maximal
irredundant quotient injects as a degreewise direct summand of the
minimal free resolution of $\oplus_{n\ge 0}\H^0(\omega_C\otimes L^n)$,
which is the dual of the minimal free resolution of
$S/I_{\varphi_{|L|}(C)}(d-g+1)$. Here and in the sequel
$S=\Sym(\H^0(L))$.  It follows that the graded betti number
$\beta_{d-2g,d-2g+2}$ of $S/I_{\varphi_{|L|}(C)}$ satisfies the
inequality
$$
\beta_{d-2g,d-2g+2}\ge \dim\Sym_{g-1}(\H^0(L\otimes\omega_C^{-1}))=
{{d-2g+1}\choose {g-1}}.$$ 
In particular $\beta_{d-2g,d-2g+2}\ge 1$,
and so $\varphi_{|L|}(C)\subset \PP^{d-g}$ doesn't satisfy property
$N_{d-2g}$, which finishes the proof of part $a)$.

Since the embedding $\varphi_{|L|}(C)\subset \PP^{d-g}$ is
arithmetically Cohen-Macaulay, its hyperplane section
$\Gamma\subset\PP^{d-g-1}$ has the same graded Betti numbers as
$\varphi_{|L|}(C)$. Thus, as in Eisenbud-Popescu [1999] or
the proof of \ref{failMRC} above, we obtain for the alternated sum
of graded betti numbers of degree $d-2g+2$ for $S/I_{\varphi_{|L|}(C)}$:
$$
\beta_{d-2g,d-2g+2}-\beta_{d-2g+1,d-2g+2}=
{{(d-g^2+g)}\over{(d-2g+2)}}{{d-g-1}\choose{g-1}}.
$$
If $3g-2\le d\le g^2-g$ the Minimal Resolution Conjecture
predicts that $\beta_{d-2g,d-2g+2}=0$, while from $a)$ we get
$\beta_{d-2g,d-2g+2}\ge 1$; hence the Minimal Resolution Conjecture 
fails for $d$ in the above range. 

For $d\ge g^2-g+1$, the expected  $\beta_{d-2g+1,d-2g+2}$ is zero while
the expected  $\beta_{d-2g,d-2g+2}$ is always larger 
than ${{d-2g+1}\choose {g-1}}$ and so the
linear complex in $a)$ doesn't account for the failure of the MRC.

For $g\ge 4$, there exist a  $g^1_{g-1}=|\cOO_C(D)|$ on $C$. The 
pairing
$$
\eta:\H^0(\cOO_C(D))\otimes \H^0(L\otimes \cOO_C(-D))\to\H^0(L) 
$$
is $1$-generic and thus $\eta$ defines a $2\times
(d-2g+2)$ matrix with linear entries in $\PP^{d-g}$ whose $2\times 2$
minors vanish in the expected codimension (see Eisenbud [1988]).
Therefore the Eagon-Northcott complex resolving the $2\times 2$-minors
of this matrix is a linear exact complex of length $d-2g+1$ which
injects as a degreewise direct summand in the top strand of the
resolution of $I_{\varphi_{|L|}(C)}$.  In particular
$\beta_{d-2g+1,d-2g+2}\ge (d-2g+1)>1$ so the Minimal Resolution
Conjecture fails for $\varphi_{|L|}(C)$ also when $d\ge g^2-g+1$. This
concludes the proof of part $b)$.
\Box

\references
\item{} M.~Auslander, D.~Buchsbaum: Codimension  and multiplicity,
{\sl Ann. of Math.} (2) {\bf 68} (1958), 625--657.
\medskip

\item{} E.~Ballico, A.V.~Geramita: The minimal free resolution of the ideal
of $s$ general points in $\PP^3$, {\sl Proceedings of the 1984
Vancouver conference in algebraic geometry}, 1--10, CMS Conf. Proc.,
{\bf 6}, Amer. Math. Soc., Providence, R.I., 1986.
\medskip

\item{} S.~Beck, M.~Kreuzer:
How to compute the canonical module of a set of points, {\sl
Algorithms in algebraic geometry and applications (Santander,
1994)}, 51--78,  Progr. Math., {\bf 143}, Birkh\"auser, Basel, 1996
\medskip

\item{} J.~Bernstein, I.N.~Gel'fand, S.I.~Gel'fand:
Algebraic bundles on $\PP^n$ and problems of linear algebra,
{\sl Funct. Anal. and its Appl.} {\bf 12} (1978), 212--214.
(Trans. from {\sl Funkz. Anal. i. Ego Priloz} {\bf 12} (1978), 66--67.)
\medskip

\item{} M.~Boij:
Artin level algebras, {\sl J. Algebra} {\bf 226} (2000), 361--374.
\medskip

\item{} I.~Dolgachev, D.~Ortland: Points sets in
projective spaces and theta functions, {\sl Ast\'erisque}
{\bf 165}, (1988).
\medskip

\item{} D.~Eisenbud, E.G.~Evans, Jr.: Generating modules efficiently:
theorems from algebraic $K$-theory, {\sl J. Algebra} {\bf 27} (1973),
278--305.
\medskip

\item{} D.~Eisenbud: Linear sections of determinantal
varieties, {\sl Amer. J. Math.} {\bf 110} (1988), 541--575.
\medskip

\item{} D.~Eisenbud, S.~Popescu: Gale Duality and Free Resolutions 
of Ideals of Points, {\sl Invent. Math.} {\bf 136} (1999), 419--449.
\medskip

\item{} D.~Eisenbud, S.~Popescu: The projective geometry of the Gale Transform,
{\sl J. Algebra.} {\bf 230} (2000), no. 1, 127--173.
\medskip

\item{} D.~Eisenbud, F.-O.~Schreyer: Sheaf Cohomology and
Free Resolutions over Exterior Algebras, preprint {\tt math.AG/0005055}.
\medskip

\item{} F.~Gaeta: Sur la distribution des degr\'es des formes
appartenant \`a la matrice de l'id\'eal homog\`ene attach\'e
\`a un groupe de $N$ points g\'en\'eriques du plan,
{\sl C.~R.~Acad.~Sci.~Paris} {\bf 233} (1951), 912--913.
\medskip

\item{} F.~Gaeta: A fully explicit resolution of the ideal
defining $N$ generic points in the plane, preprint 1995.
\medskip

\item{} A.V.~Geramita, A.M.~Lorenzini:
The Cohen-Macaulay type of $n+3$ points in $\PP^n$,
in ``The Curves Seminar at Queen's'', Vol. VI, Exp. No. F,
{\sl Queen's Papers in Pure and Appl. Math.} {\bf 83}, (1989).
\medskip

\item{} M.~Green: Koszul cohomology and the geometry of projective varieties,
{\sl J. Diff. Geom.} {\bf 19} (1984), 125--171.
\medskip

\item{} M.~Green: The Eisenbud-Koh-Stillman Conjecture on Linear Syzygies,
{\sl Invent. Math.} {\bf 136} (1999), 411--418.
\medskip

\item{} M.~Green, R.~Lazarsfeld:
Some results on the syzygies of finite sets and algebraic curves,
{\sl Compositio Math.} {\bf 67} (1988), no. 3, 301--314.
\medskip

\item{} A.~Hirschowitz, C.~Simpson: La r\'esolution minimale
d'un arrangement g\'en\'eral d'un grand nombre de points dans $\PP^n$,
{\sl Invent. Math.} {\bf 126} (1996), no. 3, 467--503.
\medskip

\item{} M.~Kreuzer: On the canonical module of a $0$-dimensional scheme,
{\sl Canadian J. of Math.} {\bf 46} (1994), no. 2, 357--379.
\medskip

\item{} F.~Lauze: preprint (1996, in preparation).
\medskip

\item{} A.M.~Lorenzini: On the Betti numbers of points in projective
space, Ph.D. thesis, Queen's University, Kingston, Ontario, 1987.
\medskip

\item{} A.M.~Lorenzini: The minimal resolution conjecture,
{\sl J. Algebra} {\bf 156} (1993), no. 1,  5--35.
\medskip

\item{} J.-P.~Serre: 
Modules projectifs et espaces fibr\'es \`a fibre vectorielle,
S\'eminaire P. Dubreil, M.-L. Dubreil-Jacotin et C. Pisot, 1957/58,
Fasc. {\bf 2}, Expos\'e 23.
\medskip

\item{} Ch.~Walter: The minimal free resolution of the homogeneous ideal of 
$s$ general points in $\PP^4$, {\sl Math. Z.} {\bf 219} (1995), no. 2, 
231--234. 
\medskip

\vskip .6cm
\widow {.1}
\vbox{\noindent Author Addresses:
\medskip
\noindent{David Eisenbud}\par
\noindent{Department of Mathematics, University of California, Berkeley,
Berkeley CA 94720}\par
\noindent{\tt de@msri.org}
\medskip
\noindent{Sorin Popescu}\par
\noindent{Department of Mathematics, SUNY at Stony Brook,
Stony Brook, NY 11794, and \par Department of Mathematics, Columbia  University,
New York, NY 10027}\par
\noindent{\tt sorin@math.sunysb.edu}\par
\medskip
\noindent{Frank-Olaf Schreyer}\par
\noindent{FB Mathematik, Universit\"at Bayreuth,
D-95440 Bayreuth, Germany}\par
\noindent{\tt schreyer@btm8x5.mat.uni-bayreuth.de}
\medskip
\noindent{Charles Walter}\par
\noindent{Laboratoire J.-A.\ Dieudonn\'e (UMR 6621 du CNRS)\par
Universit\'e de Nice -- Sophia Antipolis, 06108 Nice Cedex 02,
France}\par
\noindent{\tt walter@math.unice.fr}
}

\end